\documentclass{article}
\usepackage[utf8]{inputenc}
\usepackage[left=1in,top=1in,right=1in,bottom=1in]{geometry}
\usepackage[linktocpage=true]{hyperref}
\usepackage{setspace}
\usepackage{amssymb, amsmath, amsthm, graphicx,mathrsfs}
\usepackage{caption,cite}
\usepackage{mathtools,color}

\newtheorem{thm}{Theorem}[section]

\newtheorem{lem}[thm]{Lemma}

\newtheorem{prop}[thm]{Proposition}
\newtheorem{claim}[thm]{Claim}

\newcommand{\Fb}{\mathrm{Forb}}
\newcommand{\Hh}{\hat{H}}
\newcommand{\C}{\mathcal{C}}

\newcommand{\N}{\mathbb{N}}
\newcommand{\E}{\mathbb{E}}

\renewcommand{\P}{\mathrm{Pr}}
\renewcommand{\l}{\left}
\renewcommand{\r}{\right}
\newcommand{\ex}{\mathrm{ex}}

\usepackage[affil-it]{authblk}

\title{Random Tur\'an theorem for expansions of\\ spanning subgraphs of tight trees}
\author{Jiaxi Nie\thanks{jiaxi\_nie@fudan.edu.cn}}
\affil{Shanghai Center for Mathematical Sciences, Fudan University, Shanghai, China}
\date{\today}

\begin{document}
\maketitle

\begin{abstract}
    The $r$-expansion of a $k$-uniform hypergraph $H$, denoted by $H^{(+r)}$, is an $r$-uniform hypergraph obtained by enlarging each $k$-edge of $H$ with a set of $r-k$ vertices of degree one. The random Tur\'an number $\mathrm{ex}(G^r_{n,p},H)$ is the maximum number of edges in an $H$-free subgraph of $G^r_{n,p}$, where $G^r_{n,p}$ is the Erd\H{o}s-R\'enyi random $r$-graph with parameter $p$. In this paper, we prove an upper bound for $\mathrm{ex}(G^r_{n,p},H)$ when $H$ belongs to a large family of $r$-partite $r$-graphs: the $r$-expansion of spanning subgraphs of tight trees. This upper bound is essentially tight for at least the following two families of hypergraphs. 
    
    1. Our upper bounds are essentially tight for expansions of $K^{k-1}_{k}$, the complete $(k-1)$-graph on $k$ vertices. The proof of the lower bound makes use of a recent construction of Gowers and Janzer generalizing the famous Ruzsa-Szemer\"edi construction. In particular, when $k=3$, this answers a question of the current author, Spiro and Verstraëte concerning the random Turán number of linear triangle.

    2. Let $T$ be a tight tree such that the intersection of all edges of $T$ is empty. Simple construction shows that the upper bounds we have for expansions of $T$ are essentially tight.

    The main technical contribution of this paper is a new way to obtain balanced supersaturation results for expansions of hypergraphs: we combine two ideas, one of Mubayi-Yepremyan and another of Balogh-Narayanan-Skokan, via codegree dichotomy. We note that neither of these two ideas alone would be enough to recover results in this paper.
\end{abstract}

\section{Introduction}
An {\em $r$-uniform hypergraph}, or {\em $r$-graph} for short, is a hypergraph whose edges all have size $r$. Let $H$ be an $r$-graph, the \emph{Turán number} $\ex(n,{H})$  is the largest integer $m$ such that there exists an $H$-free $r$-graph with $n$ vertices and $m$ edges. Further, the \emph{Turán density} $\pi(H)$ is $\lim_{n\to\infty}\ex(n,H)/\binom{n}{r}$. Determining or estimating $\ex(n,{H})$ and $\pi(H)$ is a central problem in extremal combinatorics. Turán in his seminal paper\cite{turan1941} showed that $\ex(n,K_t)=(1-\frac{1}{t-1})\binom{n}{2}+\Theta(n)$, where $K_t$ is the complete graph on $t$ vertices. This result was later generalized by Erd\H{o}s, Stone and Simmonovits\cite{erdos1946structure,Erdos1966ALT} who showed that $\pi(G)=1-\frac{1}{\chi(G)-1}$, where $\chi(G)$ is the chromatic number of $G$. This result essentially solved the Tur\'an problems for all graphs with chromatic number at least three. For bipartite graphs, only sporadic results are known; see~\cite{furedi2013history} for a great survey on Tur\'an problems for bipartite graphs by F\"uredi and Simonovits. For hypergraph, it is still a major open problem to determine $\pi(K^r_k)$ for any $k>r\ge3$, where $K^r_k$ is the complete $r$-graph on $k$ vertices. See~\cite{keevash2011hypergraph} for an excellent survey by Keevash on hypergraph Turán problems. For integers $r>k\ge 2$, the {\em $r$-expansion} of a $k$-graph $H$, denoted by $H^{(+r)}$, is the $r$-graph obtained from $H$ by enlarging each $k$-edge of $H$ with a set of $r-k$ vertices of degree one. See~\cite{mubayi2016survey} for a great survey by Mubayi and Verstraëte on Turán problems for expansions. 

Let $G_{n,p}^r$ be the random subgraph of $K^r_n$ obtained by choosing each edge of $K_n^{r}$ independently with probability $p$. The random Tu\'ran number $\ex(G^r_{n,p},H)$ is the maximum number of edges in an $H$-free subgraph of $G^r_{n,p}$. The random Tu\'ran problem consists of estimating the random variable $\ex(G^r_{n,p},H)$. To this end, we say that a statement depending on $n$ holds {\em asymptotically almost surely} (abbreviated a.a.s.) if the probability that it holds tends to 1 as $n$ tends to infinity. See~\cite{rodl2013extremal} for a great survey by R\"odl and Schacht on random Tur\'an problems. For $r$-graphs that are not $r$-partite, the random Turán problems have been essentially solved independently by Conlon and Gowers~\cite{conlon2016combinatorial} and Schacht~\cite{schacht2016extremal}. To formally state their result, we define the {\em $r$-density} of $H$,
$$
m(H)=\max_{G\subset H,~e(G)\ge 2}\l\{\frac{e(G)-1}{v(G)-r}\r\}.
$$
If $p\ll n^{-1/m(H)}$, then the number of copies of $H$ in $G^r_{n,p}$ is much smaller than the number of edges in $G^r_{n,p}$, and hence $\ex(G^r_{n,p},H)=(1+o(1)e(G^r_{n,p})$. Conlon and Gowers~\cite{conlon2016combinatorial} and Schacht~\cite{schacht2016extremal} showed the following, which confirm conjectures of Kohayakawa, {\L}uczak and R\"odl~\cite{kohayakawa1997k} and R\"odl, Ruci\'nski and Schacht~\cite{rodl2007ramsey}.
\begin{thm}[\cite{conlon2016combinatorial,schacht2016extremal}]
     If $p\gg n^{-1/m(H)}$, then a.a.s.
     $$
     \ex(G^r_{n,p},H)=(\pi(H)+o(1))p\binom{n}{r}.
     $$
\end{thm}

The behaviour of $\ex(G^r_{n,p},H)$ when $H$ is $r$-partite is still a wide open problem. There are many sporadic results regarding the random Tur\'an problems for degenerate graphs and hypergraphs, see for examples~\cite{rodl2013extremal,haxell1995turan,MORRIS2016534, balogh2011number,jiang2022balanced,spiro2022random,Mubayi2023OnTR,nie2021triangle,spiro2021relative,spiro2022counting,nie2024turan}. 
In particular, the current author, Spiro and Verstra\"ete~\cite{nie2021triangle} showed the following.

\begin{thm}[\cite{nie2021triangle}]\label{theorem:C33}
    If $p\ge n^{-\frac{3}{2}+o(1)}$, then a.a.s.
    $$
    \ex(G^3_{n,p}, K^{(+3)}_3)=p^{\frac{1}{3}}n^{2+o(1)}.
    $$
\end{thm}

Moreover, for $r\ge4$, it is claimed in \cite{nie2021triangle} that if $p\ge n^{-\frac{3}{2}+o(1)}$, then a.a.s.
\begin{equation}\label{equation:triangle}
    \max\{p^{\frac{1}{2r-3}}n^{2+o(1)},pn^{r-1+o(1)}\}\le\ex(G^r_{n,p}, K^{(+r)}_3)\le   \max\{p^{\frac{3}{5}}n^{\frac{3r+3}{5}+o(1)},pn^{r-1+o(1)}\}.
\end{equation}

In this paper, we close the gap in (\ref{equation:triangle}), showing that the lower bound is essentially tight, and extend the result to expansions of $K^{k-1}_k$ for all $k\ge 3$; See Theorem~\ref{Theorem:main_clique}. In fact, we prove an upper bound for $\ex(G^r_{n,p},H)$ when $H$ belongs to a large family of $r$-partite $r$-graphs: the $r$-expansion of spanning subgraphs of tight trees. An $r$-graph $T$ is a {\em tight $r$-tree} if its edges can be ordered as $e_1,\dots,~e_t$ so that 
$$
\forall i\ge 2~\exists v\in e_i~and~1\le s\le i-1~such~that~v\not\in\cup_{j=1}^{i-1}e_j~and~e_i-v\subset e_s.
$$ 
To formally state our result, we find it more convenient to use the parameter 
$$
s(H):=\frac{1}{m(H)}=\min_{G\subset H,~e(G)\ge 2}\l\{\frac{v(G)-r}{e(G)-1}\r\}.
$$ 
Roughly speaking, $s(H)$ measures the ``spreadness'' of $H$. For integers $r\ge k\ge2$, let $T$ be a tight $k$-tree with at least two edges, let $S$ be a spanning subgraph of $T$, and let $F=S^{(+r)}$. It is not hard to check that $s(S)\ge s(T)=1$, where equality holds if and only if $S=T$. Further, $s(F)=s(S)+r-k\ge r-k+1$.
\begin{thm}\label{theorem:main}
    For integers $r\ge k\ge2$, let $\Delta=r-k$. Let $T$ be a tight $k$-tree, let $S$ be a spanning subgraph of $T$, and let $F=S^{(+r)}$. There exists a constant $C>0$ such that the following holds. Let
    $$
    p_0=Cn^{-s(F)}(\log n)^{\Delta+2+\frac{\Delta s(F)}{\Delta+1}}\text{and}~p_1=n^{-\frac{\Delta}{\Delta+1}s(F)}(\log n)^{\Delta+2+\frac{\Delta s(F)}{\Delta+1}}.
    $$
    Then a.a.s.
    $$
    \ex(G^r_{n,p},F)\le 
    \l\{
    \begin{aligned}
    &Cpn^{r-1},~~~&\text{if}~p\ge p_1;\\
    &Cp^{\frac{s(F)-\Delta-1}{s(F)}}n^{k-1}(\log n)^{\Delta+\frac{(\Delta+2)(\Delta+1)}{s(F)}},~~&\text{if}~p_0\le p\le p_1.
    \end{aligned}
    \r.
    $$
\end{thm}

In particular, Theorem~\ref{theorem:main} gives essentially tight upper bound for $\ex(G^r_{n,p},F)$ when $F$ belongs to one of the following two families of hypergraphs. 

First, for expansions of $K^{k-1}_{k}$:
\begin{thm}\label{Theorem:main_clique}
For integers $r\ge k\ge3$, let $\Delta=r-k$.
Then a.a.s
$$
    \ex(G^r_{n,p},K^{k-1(+r)}_k)=
    \l\{
    \begin{aligned}
        &\Theta(pn^{r-1}),~~&\text{if}~ p\ge n^{-\Delta-\frac{\Delta}{(\Delta+1)(k-1)}+o(1)};\\
        &p^{\frac{1}{(\Delta+1)(k-1)+1}}n^{k-1+o(1)},~~&\text{if}~n^{-\Delta-\frac{k}{k-1}+o(1)}\le p\le n^{-\Delta-\frac{\Delta}{(\Delta+1)(k-1)}+o(1)}.
    \end{aligned}
    \r.
    $$
\end{thm} 

The lower bound construction of Theorem~\ref{Theorem:main_clique} makes use of a recent construction of Gowers and Janzer generalizing the famous Ruzsa-Szemer\"edi construction.

\begin{figure}[h]
    \centering
    \includegraphics[scale=0.4]{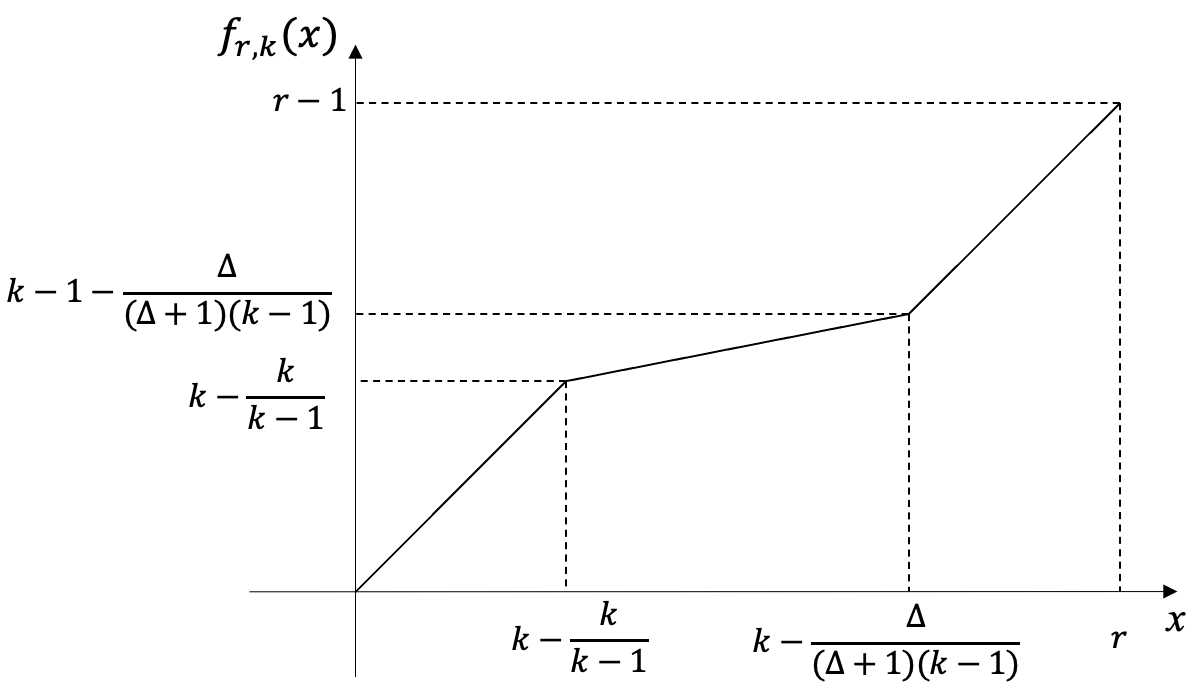}
    \caption{The behaviour of $\ex(G^r_{n,p}, K^{k-1(+r)}_k)$ with respect to $p$.}
    \label{figure:clique}
\end{figure}

Let $0< x\le r$, and let $p=n^{-r+x}$. Define
\begin{equation*}
f_{r,k}(x)=\lim_{n\rightarrow\infty}\log_n\ex(G^{r}_{n,p}, K^{k-1(+r)}_k).
\end{equation*}
Figure~\ref{figure:clique} roughly describes the behaviour of $\ex(G^r_{n,p}, K^{k-1(+r)}_k)$ with respect to $p$.

It is very interesting that the middle range here is not ``flat'', since, as far as we know, it is always the case that there is a ``flat'' middle range for other degenerate random Turán problems that have been essentially solved. A recent paper of the current author and Spiro~\cite{nie2023sidorenko} showed that such phenomenon is closely related to the fact that $K^{k-1(+r)}_k$ does not have the so-called Sidorenko property; see \cite{nie2023sidorenko} for more explanations.

The next result is another example where the middle range is ``flat''. For expansions of tight $k$-trees whose edges do not contain a common vertex, we prove the following.
\begin{thm}\label{theorem:main_tighttree}
    For $r\ge k\ge2$, let $\Delta=r-k$ and let $T$ be a tight $k$-tree such that $\cap_{e\in E(T)}e=\emptyset$. Then a.a.s.
    $$
    \ex(G^r_{n,p}, T^{(+r)})=\l\{
    \begin{aligned}
        &\Theta(pn^{r-1}),~~&\text{if}~p\ge n^{-\Delta+o(1)};\\
        &n^{k-1+o(1)},~~&\text{if}~n^{-\Delta-1+o(1)}\le p\le n^{-\Delta+o(1)}.
    \end{aligned}
    \r.
    $$
\end{thm}

\begin{figure}[h]
    \centering
    \includegraphics[scale=0.4]{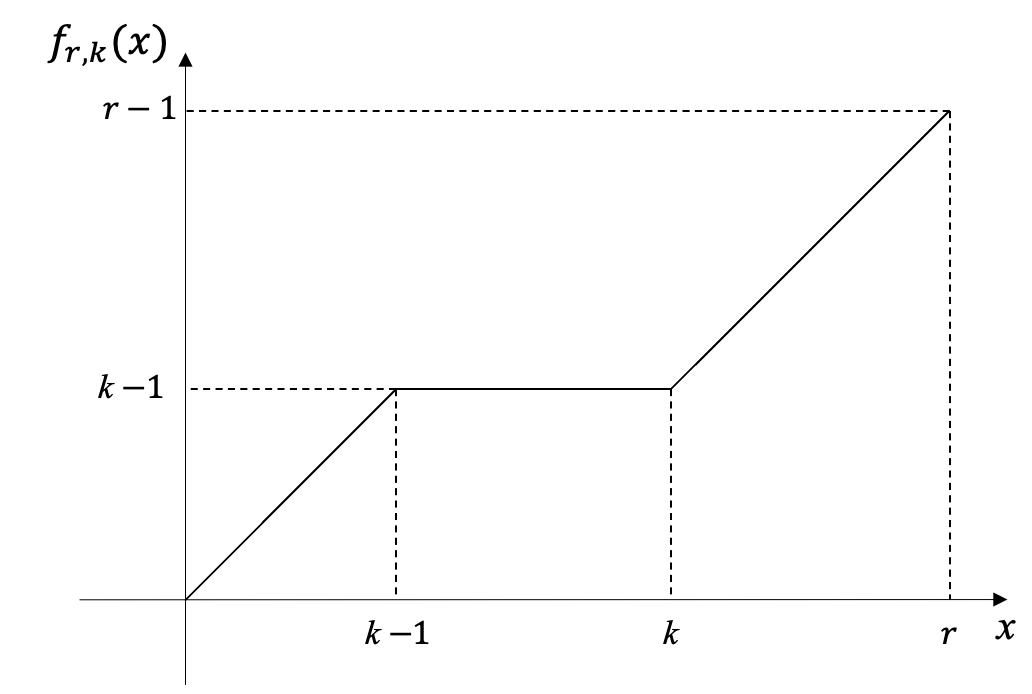}
    \caption{The behaviour of $\ex(G^r_{n,p}, T^{(+r)})$ with respect to $p$.}
    \label{figure:tree}
\end{figure}

Let $0< x\le r$, and let $p=n^{-r+x}$. Define
\begin{equation*}
f_{r,k}(x)=\lim_{n\rightarrow\infty}\log_n\ex(G^{r}_{n,p}, T^{(+r)}).
\end{equation*}
Figure~\ref{figure:tree} roughly describes the behaviour of $\ex(G^r_{n,p}, T^{(+r)})$ with respect to $p$.

We would like to note that Theorem~\ref{theorem:main} is not always optimal. Let $C^r_\ell:=C^{(+r)}_{\ell}$ denote the $r$-uniform linear cycles of length $\ell$.  Theorem~\ref{theorem:main} can also be applied on linear cycles, but it does not give the optimal upper bound. Recently, the current author~\cite{nie2024turan} and, independently, Mubayi and Yepremyan~\cite{Mubayi2023OnTR} essentially solved the random Tur\'an problem for $C^r_{2\ell}$ when $r\ge 4$.

\begin{thm}[\cite{nie2024turan,Mubayi2023OnTR}]\label{theorem:main_ge4}
    For every $r\ge 4$ and $\ell\ge2$, if $p\ge n^{-r+2+\frac{1}{2\ell-1}+o(1)}$, then a.a.s.
    $$
    \ex(G^r_{n,p}, C^r_{2\ell})\le pn^{r-1+o(1)}.
    $$
\end{thm}

When $r=3$, the current author~\cite{nie2024turan} showed the following upper bound, which is essentially tight for all but a small range of $p$; see~\cite{nie2024turan} for more.

\begin{thm}[\cite{nie2024turan}]\label{theorem:main_3}
    For every $\ell\ge2$, there exists $c>0$ such that the following holds. Let
    $$
    p_0=cn^{-\frac{4\ell-3}{4\ell-2}}(\log n)^{\frac{4\ell^2+\ell-4}{2\ell-2}},~p_1=n^{-\frac{(\ell-1)(4\ell-3)}{4\ell^2-5\ell+2}}(\log n)^{3+\frac{2\ell-2}{4\ell^2-5\ell+2}}.
    $$
    Then a.a.s.
    $$
    \ex(G^3_{n,p}, C^3_{2\ell})\le 
    \l\{
    \begin{aligned}
        &cp^{\frac{2(\ell-1))}{\ell(4\ell-3)}}n^{1+\frac{1}{\ell}}(\log n)^{\frac{2(2\ell-1)^2}{\ell(4\ell-3)}},~~&\text{if}~p_0\le p<p_1;\\
        &cpn^2,~~&\text{if}~p\ge p_1.
    \end{aligned}
    \r.
    $$
\end{thm}

For linear odd cycle, Theorem~\ref{theorem:main} gives the following upper bound.

\begin{thm}\label{theorem:main_oddcycle}
    For $r\ge 3$ and $\ell\ge 1$, there exists constant $c$ such that a.a.s.
    $$
    \ex(G^r_{n,p},C^r_{2\ell+1})\le
    \l\{
    \begin{aligned}
        &cpn^{r-1},~~&\text{if}~ p\ge n^{-r+2+\frac{2\ell+r-3}{2\ell(r-2)}+o(1)};\\
        &p^{\frac{2\ell-1}{2\ell(r-1)-1}}n^{2+o(1)},~~&\text{if}~n^{-r+1+\frac{1}{2\ell}+o(1)}\le p\le n^{-r+2+\frac{2\ell+r-3}{2\ell(r-2)}+o(1)}.
    \end{aligned}
    \r.
    $$
    Moreover, when $\ell=1$, this upper bound is tight.
\end{thm}

The lower bound for $\ell=1$ comes from Theorem~\ref{Theorem:main_clique} when $k=3$. When $\ell\ge2$, we believe that this upper bound is not tight, and we think it is an interesting problem to determine the magnitude of $\ex(G^r_{n,p},C^r_{2\ell+1})$.

Theorem~\ref{theorem:main} relies on a ``balanced supersaturation'' result for expansions of spanning subgraphs of tight trees. Here ``balanced supersaturation'' roughly means that, given an $r$-graph $H$, if the number of edges in an $r$-graph $G$ on $n$ vertices is much larger than $\ex(n,H)$, then we can find a collection $\C$ of copies of $H$ in $G$ with ``balanced distribution'', i.e. you cannot find a set $\sigma$ of edges of $G$ such that the number of copies of $H$ in $\C$ containing $\sigma$ is unusually large. A good balanced supersaturation result for $H$ would usually imply a good upper bound for the random Tur\'an number of $H$ by using the container method developed independently by Balogh, Morris and Samotij~\cite{balogh2015independent} and Saxton and Thomassen~\cite{saxton2015hypergraph}. 

Let $\Fb(H,n)$ be the number of $H$-free hypergraphs on $n$ vertices. Determining $\Fb(H,n)$ is another type of problem that relies heavily on balanced supersaturation results. This type of problem has been studied extensively for $r$-graphs that are not $r$-partite~\cite{erdos1986asymptotic,nagle2001asymptotic,nagle2006extremal}. Mubayi and Wang~\cite{mubayi2019thenumber} initiated the study for $r$-partite $r$-graphs. In~\cite{mubayi2019thenumber}, they determined the asymptotics of $\Fb(C^r_k,n)$ for even $k$ and $r=3$, and conjectured the asymptotics for all $r$ and $k$. Their conjecture was later confirmed by Balogh, Skokan and Narayanan~\cite{balogh2019number}. Ferber, Mckinley and Samotij~\cite{ferber2020supersaturated} proved similar results for a larger family of hypergraphs which includes linear cycles. The results of~\cite{balogh2019number} and \cite{ferber2020supersaturated} both rely on balanced supersaturation. More recently, Jiang and Longbrake~\cite{jiang2022balanced} prove the balanced supersaturation result in~\cite{ferber2020supersaturated} in a more explicit form. Our main technical contribution in this paper is that we find a new way to obtain balanced supersaturation results for expansions of hypergraphs, which significantly improve the balanced supersaturation results in~\cite{balogh2019number,ferber2020supersaturated,jiang2022balanced}.

\subsection{Proof outline}
The main result in this paper is an upper bound for $\ex(G^r_{n,p},H)$ when $H$ is an $r$-expansion of a spanning subgraph of a tight tree, which mainly comes from proving a novel balanced supersaturation theorem. Obtaining an upper bound of random Tur\'an number from a balanced supersaturation result is by far a routine application of the container method developed independently by Balogh, Morris and Samotij~\cite{balogh2015independent} and Saxton and Thomassen~\cite{saxton2015hypergraph}. Hence we will only outline the proof of the new balanced supersaturation theorem. 

Let $T$ be a tight $k$-tree, let $S$ be a spanning subgraph of $T$, and let $F=S^{(+r)}$ be the $r$-expansion of $S$. Recall that ``balanced supersaturation'' roughly means that if the number of edges in an $r$-graph $H$ on $n$ vertices is much larger than $\ex(n,F)$, then we can find a collection $\C$ of copies of $F$ in $H$ with ``balanced distribution''. To construct a balanced collection of copies of $F$, we will adopt the following strategy using codegree dichotomy and induction on $r$. A {\em $k$-shadow} of $H$ is a set $\sigma$ of $k$ vertices such that $\sigma\subset e$ for some edge $e$ of $H$. The {\em codegree} of $\sigma$ is the number of edges of $H$ containing $\sigma$.

\begin{figure}[h]
    \centering
    \includegraphics[scale=0.5]{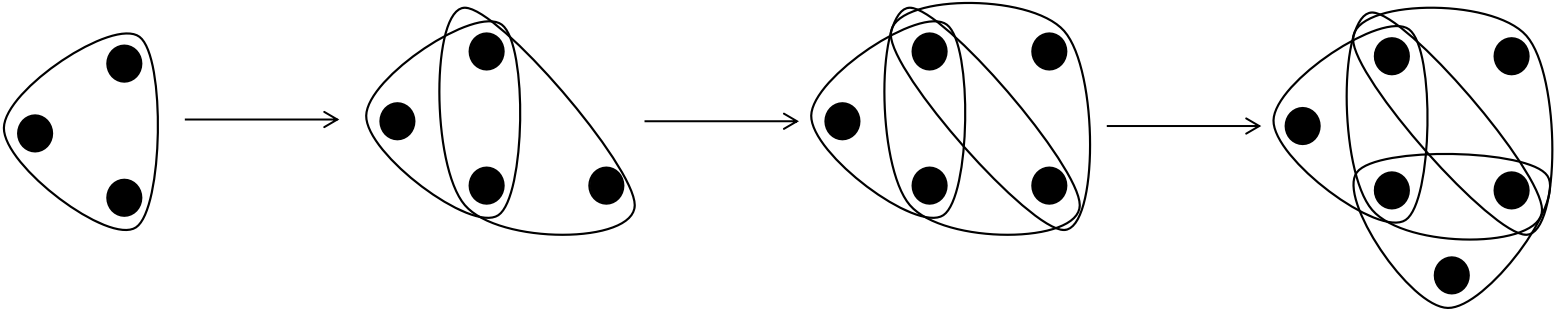}
    \caption{Constructing $C^3_3$ using ``greedy expansion''}
    \label{fig:greedy}
\end{figure}

For the base case $r=k$, given an $r$-graph $H$ with $tn^{r-1}$ edges where $t$ is sufficiently large, we can first obtain a subgraph $H'$ of $H$ with $\Omega(tn^{r-1})$ edges such that every $(r-1)$-shadow of $H'$ has codegree at least $t$ by deleting unwanted edges. Then we can construct a balanced collection of copies of $F$ in $H'$ using ``greedy expansion'' (see Figure~\ref{fig:greedy}): we greedily embed the edges of the tight tree $T$ into $H'$; since every $(r-1)$-shadow of $H'$ has codegree at least $t$, there are $\Omega(t)$ ways to embed each new edge. Once we obtain a copy of $T$, we automatically obtain a copy of $F$ because $F$ is a spanning subgraph of $T$. With some additional techniques, we can ensure that this collection of copies of $F$ generated greedily is evenly distributed in $H'$. This ``greedy expansion'' method is inspired by Balogh, Narayanan and Skokan~\cite{balogh2019number}. 

For $r>k$, suppose we already have a balanced supersaturation theorem for $F'=S^{(+(r-1))}$ and we want to construct a balanced collection of copies of $F=S^{(+r)}$. Given an $r$-graph $H$ with many edges, we may and will assume that it is $r$-partite with respect to the partition $V(H)=V_1\cup\cdots\cup V_r$. The proof splits into two cases according to how large is the ``typical'' codegree of the $(r-1)$-shadows of $H$. More precisely, $H$ can be partitioned into subgraphs $\Hh_L$ and $\Hh_{\tau, a}$, $\tau\in[r]^{r-1}$, $a\ge 0$, such that: 1. every $(r-1)$-shadows of $\Hh_L$ is ``large'', and 2. every $(r-1)$-shadow in $\cup_{i\in\tau}V_i$ of $\Hh_{\tau,a}$ has codegree at least $2^a$ and at most $2^{a+1}$ where $2^a$ is ``small''. The idea of partitioning $H$ into subgraphs with ``almost regular'' codegrees can be found in~\cite{mubayi2020random}. The proof splits into two cases as following.

Case 1: if $\Hh_L$ has more than $e(H)/(\log n)^{r-1}$ edges, then again we can construct a balanced collection of $F$ using ``greedy expansion''. Because the codegree of every $(r-1)$-shadow is ``large'', we can obtain many copies of $F$ in this way. 

Case 2: if $e(\Hh_L)$ is smaller than $e(H)/(\log n)^{r-1}$ edges, then we can find $(r-1)$ parts, say $V_1,\dots,V_{r-1}$, and a subgraph $\Hh$ of $H$ with $e(\Hh)\ge e(H)/\log n$ such that all $(r-1)$-shadows of $\Hh$ in $V_1\cup\cdots\cup V_{r-1}$ have codegree at least $D$ and at most $2D$ where $D$ is ``small''. We now consider the $(r-1)$-graph $G$ on $V_1\cup\dots\cup V_{r-1}$ whose edge set consists of all $(r-1)$-shadows of $\Hh$ in $V_1\cup\cdots\cup V_{r-1}$. Then $e(G)\ge e(\Hh)/(2D)$. Since $D$ is ``small'', $e(G)$ must be ``large''. Applying the balanced supersaturation theorem for $F'$ on $G$ gives a balanced collection $\C'$ of copies of $F'$ in $G$. Using the fact that the codegrees are ``almost regular'', we can extend each copy of $F'$ in $\C'$ into copies of $F$ in $H$, which gives a balanced collection $\C$ of copies of $F$ in $H$. This idea of first using the supersaturation in the shadow and then extending can be traced back to Mubayi and Yepremyan~\cite{mubayi2020random}.

By selecting the best cutoff that defines ``large'' and ``small'' codegrees, we obtain the balanced supersaturation theorems in this paper. To summarize, given an $r$-graph $H$ with many edges, if the typical codegrees of its $(r-1)$-shadows are ``large'', then we find copies of $F$ by ``greedy expansion''; if the typical codegrees are ``small'', then we find copies of $F'$ in its $(r-1)$-shadows and then expand them into copies of $F$. This idea of codegree dichotomy is the main novelty of this paper.

\subsection{Notations and structure}
Let $H=(V,E)$ be an $r$-graph, we denote $|V|$ by $v(H)$ and $|E|$ by $e(H)$. Given $E'\subset E$, $H-E'$ is the $r$-graph $(V,E\setminus E')$. A {\em $k$-shadow} of $H$ is a set $\sigma$ of $k$ vertices such that $\sigma\subset e$ for some edge $e$ of $H$. Given a set $\sigma$ of vertices of $H$, let $N_H(\sigma)$ be the set of edges of $H$ containing $\sigma$ and let $d_{H}(\sigma)=|N_H(\sigma)|$. For $1\le k\le r$, let $\Delta_k(H)$ denote the maximum $d_H(\sigma)$ with $|\sigma|=k$. Given a collection $\C$ of copies of $C^r_{2\ell}$ in an $r$-graph $H$, we can view $\C$ as a hypergraph on $E(H)$ whose edges are copies of $C^r_{2\ell}$ in $\C$. Then $\Delta_k(\C)$ means the maximum number of copies of $C^r_{2\ell}$ in $\C$ containing $k$ common edges. Given two positive functions $f(n)$ and $g(n)$ on $\N$, we write $f(n)\le O(g(n))$ if there exists a constant $c>0$ such that $f(n)= c\cdot g(n)$ for all $n$, and write $f(n)= \Omega(g(n))$ if $g\le O(f)$. If $f= O(g)$ and $g= O(f)$, then we write $f=\Theta(g)$. We write $f(n)=o(g(n))$ if $f(n)/g(n)\to 0$ as $n\to\infty$, and write $f(n)=\omega(g(n))$ if $g(n)=o(f(n))$. 

This paper is strutured as follows. In Section~\ref{section:bss}, we prove a balanced supersaturation theorem (Theorem~\ref{theorem:bss_main}) for expansions of spanning subgraphs of tight trees. In Section~\ref{section:container}, we use the hypergraph container method, together with Theorem~\ref{theorem:bss_main}, to prove a container theorem for expansions of spanning subgraphs of tight trees. In Section~\ref{section:main}, we prove Theorem~\ref{theorem:main}, the random Turán upper bound for expansions of spanning subgraphs of tight trees. In Section~\ref{section:construction}, we introduce two simple constructions for all hypergraphs and a nontrivial construction for expansions of $K^{k-1}_k$. These constructions, together with Theorem~\ref{theorem:main}, imply Theorem~\ref{Theorem:main_clique} and Theorem~\ref{theorem:main_tighttree}.

\section{Balanced Supersaturation}\label{section:bss}
We first prove a balanced supersaturation result (Lemma~\ref{lemma:bss_sub_tree}) for spanning subgraphs of tight trees using the idea of ``greedy expansion'', which is inspired by Theorem 3.1 in~\cite{balogh2019number} of Balogh, Narayanan, and Skokan. 

We need the following folklore proposition.
\begin{prop}\label{prop:regular}
    For any integers $n>t>0$, there exists a graph on $n$ vertices such that each vertex has degree $t$ or $t-1$.
\end{prop}

\begin{lem}\label{lemma:bss_sub_tree}
    Let $T$ be a tight $r$-tree and let $F$ be a spanning subgraph of $T$. There exist a constant $C>0$ such that the following holds for all sufficiently large $n$ and every $t$ with $2v(T)\le t\le \binom{n}{r}/n^{r-1}$. Let $H$ be an $r$-graph on $n$ vertices. If every $(r-1)$-shadow $\sigma$ of $H$ has $d_H(\sigma)> t$, then there exists a collection $\C$ of copies of $F$ in $H$ satisfying
    \begin{enumerate}
        \item[(a)] $|\C|\ge C\cdot e(H)t^{v(F)-r}$;
        \item[(b)] $\Delta_j(\C)\le C^{-1}t^{v(F)-r-(j-1)s(F)}$,~~$1\le j\le e(F)$.
    \end{enumerate}
\end{lem}

\begin{proof}
    By Proposition~\ref{prop:regular}, for each $(r-1)$ shadow $\sigma$ of $H$ we can define a graph $\Gamma(\sigma)$ on $N_H(\sigma)$ such that $\forall e\in N_H(\sigma)$ we have $d_{\Gamma(\sigma)}(e)=t$ or $t-1$. Let $\ell=e(T)$. Since $T$ is a tight tree, its edges can be ordered as $e_1,\dots,e_{\ell}$ such that 
    \begin{equation}\label{equation:1_bss_sub_tree}
    \forall 2\le i\le \ell,~\exists v_i\in e_i~\text{and}~1\le s_i\le i-1 \text{~such that~} v_i\not\in\cup_{j=1}^{i-1}e_j~\text{and}~e_i-v_i\subseteq e_{s_i}.    
    \end{equation}

    We now describe an algorithm which constructs copies of $T$ in $H$; it involves specifying $\ell$ edges $e_1,\dots,e_{\ell}$.
    \begin{enumerate}
        \item[(i)] We start by choosing $e_1$ in $E(H)$;
        \item[(ii)] For $2\le i\le \ell$, suppose we have already specified $e_1,\dots,e_{i-1}$. By definition, $e_i-v_i$ is contained in a specified edge $e_{s_i}$, where $1\le s_i\le i-1$. We then specify the vertex $v_i$ such that $v_i$ is distinct from all previously specified vertices and that $e_i$ and $e_{s_i}$ form an edge in $\Gamma(e_i\cap e_{s_i})$. (See Figure~\ref{fig:greedy})
    \end{enumerate}

    Let $\C$ be the collection of all copies of $F$ contained in copies of $T$ generated by the algorithm above. Note that for step (i) of the algorithm, the number of ways to do it is $e(H)$; for step (ii), the number of ways to specify each new vertex is at least $t/2$ (since $t\ge 2v(T)$) and at most $t$ (by the definition of $\Gamma(e_i\cap e_{s_i})$).
    \begin{claim}
    $$
    |\C|\ge\Omega(e(H)t^{v(F)-r})
    $$
    \end{claim}
    \begin{proof}
        Clearly, the number of copies of $T$ generated by the algorithm is $\Omega(e(H)t^{v(T)-r})$. Since each copies of $T$ contains a constant number of copies of $F$ and each copies of $F$ is contained in at most a constant number of $T$, we have $|\C|\ge\Omega(e(H)t^{v(F)-r})$.
    \end{proof}
    \begin{claim}
    For $1\le j\le \ell$,
    $$
    \Delta_j(\C)\le O(t^{v(F)-r-(j-1)s(F)}).
    $$
    \end{claim}
    \begin{proof}
    Let $J$ be a set of $j$ edges of $H$. We want to bound from above the number of ways to specify a copy of $F$ in $\C$ containing $J$. Clearly, it suffices to bound from above the number of ways to specify a copy of $T$ generated by the algorithm containing $J$. First, we fix which edges of $J$ correspond to which edges of $T$; this can be done in $O(1)$ ways. Note that the number of vertices fixed by $J$ is at least
    $$
    \min_{F'\subset F,~e(F')=j}v(F')= r+(j-1)\frac{\min_{F'\subset F,~e(F')=j}v(F')-r}{j-1}\ge r+(j-1)s(F).
    $$ 
    Thus the number of un-specified vertices is at most $v(F)-r-(j-1)s(F)$

    Consider a tree $\mathcal{T}$ whose vertices are the edges of a copy of $T$ generated by the algorithm and, for $2\le i\le \ell$, $e_i$ and $e_{s_i}$ form an edge in $\mathcal{T}$ (See~(\ref{equation:1_bss_sub_tree}) to recall the definition of $s_i$). According to the algorithm, if two edges $e$ and $f$ are adjacent in $\mathcal{T}$, then they are adjacent in $\Gamma(e\cap f)$. Thus by the definition of $\Gamma(e\cap f)$, if $e$ is fixed, then there are at most $t$ ways to specify $f$.

    Clearly, if $T$ has not been completely specified, then there exists an un-specified edge $f$ and a specified edge $e$ such that $f$ and $e$ form an edge in $\mathcal{T}$. As discussed above, there are at most $t$ ways to specify $f$. Note that every time we specify a new edge, the number of specified vertices decreases by one. Hence the copy of $T$ will be completely specified after specifying at most $v(F)-r-(j-1)s(F)$ such edges $f$. Therefore we have $\Delta_j(\C)\le O(t^{v(F)-r-(j-1)s(F)})$.
    \end{proof}
    \end{proof}

In the rest of this section, we will use Lemma~\ref{lemma:bss_sub_tree} as a base case, together with the idea of codegree dichotomy, to prove the following stronger balanced supersaturation result for expansions of spanning subgraph of tight tree.

\begin{thm}\label{theorem:bss_main}
For $r\ge k\ge2$, let $\Delta=r-k$. Let $T$ be a tight $k$-tree, let $S$ be a spanning subgraph of $T$, and let $F=S^{(+r)}$. There exist constants {$K,C>0$} such that the following holds for all sufficiently large $n$ and every $t$ with $K\le t\le \binom{n}{r}/n^{r-1}$. Let $H$ be an $r$-graph on $n$ vertices with $tn^{r-1}$ edges and let 
$$
A=t^{\frac{1}{\Delta+1}}\l(\frac{n}{\log n}\r)^{\frac{\Delta}{\Delta+1}}.
$$
Then there exists a collection $\C$ of copies of $F$ in $H$ such that
$$
\Delta_j(\C)\le\frac{C(\log n)^\Delta|\C|}{tn^{r-1}}A^{-s(F)(j-1)}.
$$
\end{thm}

\begin{proof}
    We prove this theorem by induction on $\Delta$. When $\Delta=0$, $F=S$ is a spanning subgraph of $T$. We want to find a subgraph $H'$ of $H$ such that every $(r-1)$ shadow of $H'$ has large codegree. To this end, we iteratively delete edges containing an $(r-1)$-shadow with codegree at most $t/2$. More precisely, we run the following algorithm:
    \begin{itemize}
        \item[(i)] Let $H_0=H$. For $i\ge0$, given $H_i$, if $H_i$ contains an $(r-1)$-shadow $\sigma$ such that $d_{H_i}(\sigma)\le t/2$, we arbitrarily fix such a $\sigma$ and let $H_{i+1}=H_i-N_{H_i}(\sigma)$. 
        \item[(ii)] If every $(r-1)$-shadow $\sigma$ of $H_i$ has $d_{H_i}(\sigma)>t/2$, then we let $H'=H_i$.
    \end{itemize}
    In step (i), we delete at most $t/2$ edges for each $(r-1)$ shadow. Hence, the number of edges deleted in step (i) is at most $\binom{n}{r-1}t/2\le tn^{r-1}/2$. Hence, we have $e(H')\ge e(H)-tn^{r-1}/2\ge e(H)/2$. Note that every $(r-1)$-shadow $\sigma$ of $H'$ has $d_{H'}(\sigma)>t/2$ and that $t\ge K$ can be sufficiently large, we can apply Lemma~\ref{lemma:bss_sub_tree} on $H'$ to obtain a collection $\C$ of copies of $F$ in $H'$ such that
    $$
    |\C|\ge\Omega\l((tn^{r-1})t^{v(F)-r}\r),
    $$
    and
    $$
    \Delta_j(\C)\le O\l(t^{v(F)-r-(j-1)s(F)}\r)\le O\l(\frac{|\C|}{tn^{r-1}}t^{-s(F)(j-1)}\r),~\forall 1\le j\le e(F).
    $$
    This completes the proof for $\Delta=0$.
    
    When $\Delta\ge 1$, suppose we have proved this theorem for $\Delta-1$. We pick an $r$ partition of $V(H)=V_1\cup\dots\cup V_{r}$ uniformly at random and let $H'$ be an $r$-partite subgraph of $H$ induced by the partition. Clearly, the expectation of $e(H')$ is $\frac{r\,!}{r^r}tn^{r-1}$. So we can fix an $H'$ with at least $\frac{r\,!}{r^r}tn^{r-1}$ edges. 
    
    Next, we describe an algorithm that partitions $H'$ into subgraphs $\Hh_L$ and $\Hh_{\tau, a}$, where $\tau\in \binom{[r]}{r-1}$, $0\le a< \log_2 A$.
    
    \begin{itemize}
        \item[(i)] Let $H_0=H$. For $i\ge 0$, given $H_i$, if there is an $(r-1)$-shadow $\sigma$ of $H_{i}$ such that $d_{H_{i}}(\sigma)\le A$, then we arbitrarily fix such a $\sigma$ and let $H_{i+1}=H_{i}-N_{H_{i}}(\sigma)$. There is a unique pair of $(\tau,a)$ such that $\sigma$ contain a vertex of $V_j$ for each $j\in\tau$ and $2^a\le d_{H_{i}}(\sigma)<2^{a+1}$. We put all the elements in $N_{H_{i-1}}(\sigma)$ into the edge set of $\Hh_{\tau, a}$. 
        \item[(ii)] If all $(r-1)$-shadows $\sigma$ of $H_{i}$ have $d_{H_{i}}(\sigma)> A$, then we let $\Hh_L=H_{i}$.
    \end{itemize}
    The proof splits into two cases.\\
    
    \noindent\textbf{Case 1: $e(\Hh_L)\ge e(H')/(\log n)^{\Delta}$.}
    
    Since every $(r-1)$-shadow $\sigma$ of $\Hh_L$ has $d_{\Hh_L}(\sigma)>A$ and $A$ is sufficiently large, we can apply Lemma~\ref{lemma:bss_sub_tree} on $\Hh_L$ to obtain a collection $\C$ of copies of $F$ in $\Hh_L$ such that
    $$
    |\C|\ge\Omega\l(\frac{tn^{r-1}}{(\log n)^{\Delta}}A^{v(F)-r}\r),
    $$
    and
    $$
    \Delta_j(\C)\le O\l(A^{v(F)-r-(j-1)s(F)}\r)\le O\l(\frac{(\log n)^{\Delta}|S|}{tn^{r-1}}A^{-s(F)(j-1)}\r),~\forall 1\le j\le e(F).
    $$
    
    \noindent\textbf{Case 2: $e(\Hh_L)<e(H')/(\log n)^{\Delta}$.}

    Note that the number of edges in $\Hh_{\tau,a}$ is small when $a$ is small. We consider the number of edges in all $\Hh_{\tau,a}$ with $2^a\le \frac{r\,!}{4r^r}t$. By step (i) of the algorithm, we know that each $(r-1)$-shadow contributes at most $\frac{r\,!}{4r^r}t$ to the following sum:
    $$
    \sum_{\tau\in \binom{[r]}{r-1}}\sum_{2^a\le \frac{r\,!}{4r^r}t} e(\Hh_{\tau,a})\le \binom{n}{r-1}\frac{r\,!}{2r^r}t\le \frac{r}{2r^r} tn^{r-1}\le e(H')/2.
    $$
    Thus the number of edges in all $\Hh_{\tau,a}$ with $2^a> \frac{r\,!}{4r^r}t$ is large, that is,
    $$
    \sum_{\tau\in \binom{[r]}{r-1}}\sum_{2^a>\frac{r\,!}{4r^r}t} e(\Hh_{\tau,a})\ge \Omega(tn^{r-1}).
    $$
    
    By the Pigeonhole Principle there exist $\tau$ and $a$ with $2^a> \frac{r\,!}{4r^r}t$ such that $e(\Hh_{\tau,a})\ge\Omega(tn^{r-1}/\log n)$. Fix such a pair of $\tau$ and $a$, let $\Hh=\Hh_{\tau,a}$ and let $D=2^a$. Since $D\ge \frac{r\,!}{4r^r}t$, $D$ would be sufficiently large if $t\ge K$ is sufficiently large. Without lost of generality, we may and will let $\tau=\{1,\dots,~r-1\}$. Let $G$ be the $(r-1)$-graph on $V_1\cup\dots\cup V_{r-1}$ whose edge set consists of all $(r-1)$-shadows of $\Hh$ in $V_1\cup\dots\cup V_{r-1}$. Then $e(G)\ge e(\Hh)/(2D)\ge\Omega(tn^{r-1}/(\log n\cdot D))$. Let
    \begin{equation}\label{equation:1_bss_main}
    t'=\frac{e(G)}{n^{r-2}}\ge\Omega\l(\frac{tn}{\log n\cdot D}\r),    
    \end{equation}
    let $F'=S^{(+(r-1))}$ and let
    \begin{equation}\label{equation:2_bss_main}
    A'=t'^{\frac{1}{\Delta}}\l(\frac{n}{\log n}\r)^{\frac{\Delta-1}{\Delta}}\ge \Omega\l(A^{\frac{\Delta+1}{\Delta}}D^{-\frac{1}{\Delta}}\r).
    \end{equation}
    The above inequality is obtained by using (\ref{equation:1_bss_main}) and the definition of $A$. By the inductive hypothesis together with (\ref{equation:1_bss_main}) and (\ref{equation:2_bss_main}), there exists a collection $\C'$ of copies of $F'$ in $G$ satisfying, $\forall 1\le j\le e(F)$,
    \begin{equation}\label{equation:3_bss_main}
    \Delta_j(\C')\le O\l(\frac{(\log n)^{\Delta-1}|\C'|}{t'n^{r-2}}A'^{-s(F')(j-1)}\r)\le O\l(\frac{(\log n)^{\Delta}|\C'|}{tn^{r-1}}D^{1+\frac{(s(F)-1)(j-1)}{\Delta}}A^{-\frac{\Delta+1}{\Delta}(s(F)-1)(j-1)}\r).    
    \end{equation}
     Note that each $(r-1)$-edge of $G$ is contained in at least $D$ $r$-edges of $\Hh$. Given that $t\ge K$ is sufficiently large (hence $D$ is sufficiently large), we can greedily extend a copy of $F'$ in $G$ into a copy of $F$ in $\Hh$ in $\Omega(D^{e(F)})$ ways. More formally, let $\C$ be the collection of copies of $F$ in $\Hh$ with the following property: the $(r-1)$-shadows of the copy of $F$ in $V_1\cup\dots\cup V_{r-1}$ form a copy of $F'$ in $\C'$. Then we have
    \begin{equation}\label{equation:4_bss_main}
    |\C|\ge\Omega\l(|\C'|D^{e(F)}\r)        
    \end{equation}
    For any $1\le j\le e(F)$, given a $j$-tuple of $r$-edges of $H'$, if they are contained in some copies of $F$ in $\C$, then the $j$ $(r-1)$-shadows of all these $j$ $r$-edges are distinct. By definition, the number of copies of $F'$ in $\C'$ containing these $j$ $(r-1)$-shadows is at most $\Delta_j(\C')$. Since each $(r-1)$-shadows is contained in at most $2D$ $r$-edges of $H'$, we have
    \begin{equation}\label{equation:5_bss_main}
            \Delta_j(\C)\le O\l(\Delta_j(\C')D^{e(F)-j}\r).
    \end{equation}

    Combining (\ref{equation:3_bss_main}), (\ref{equation:4_bss_main}), (\ref{equation:5_bss_main}), we conclude that, $\forall 1\le j\le e(F)$,
    $$
    \begin{aligned}
        \Delta_j(\C)&\le O\l(\frac{(\log n)^{\Delta}|\C'|}{tn^{r-1}}D^{e(F)-j+1+\frac{(s(F)-1)(j-1)}{\Delta}}A^{-\frac{\Delta+1}{\Delta}(s(F)-1)(j-1)}\r)\\
        &\le O\l(\frac{(\log n)^{\Delta}|\C|}{tn^{r-1}}A^{-s(F)(j-1)}\l(\frac{D}{A}\r)^{\frac{(s(F)-\Delta-1)(j-1)}{\Delta}}\r).
    \end{aligned}
    $$
    Further, since $D\le A$ and $s(F)\ge \Delta+1$,
    We have
    $$
    \Delta_j(\C)\le O\l(\frac{(\log n)^{\Delta}|\C|}{tn^{r-1}}A^{-s(F)(j-1)}\r)
    $$
\end{proof}

\section{Container Theorems}\label{section:container}
We make use of the hypergraph container method developed independently by Balogh, Morris and Samotij~\cite{balogh2015independent} and Saxton and Thomassen~\cite{saxton2015hypergraph}. More precisely, we use the following simplified version of the hypergraph container theorem of Balogh, Morris and Samotij~\cite{balogh2018method}.
\begin{thm}[\cite{balogh2018method}]\label{theorem:container}
For every $r\ge 2$, there exists constant $\epsilon>0$ such that the following holds.
Let $H$ be an $r$-graph on $n$ vertices such that
    \begin{equation}\label{equation:container}
        \triangle_j(H)\leq \l(\frac{B}{n}\r)^{j-1}\frac{|H|}{L},~~\forall 1\le j\le r, 
    \end{equation}
    for some integers $B, L>0$.
    Then there exists a collection $\mathcal{C}$ of at most
    $$
    \exp\l(\frac{\log(\frac{n}{B})B}{\epsilon}\r)
    $$
    subsets of $V(H)$ such that:
    \begin{enumerate}
        \item[(a)] for every independent set $I$ of $H$, there exists $C\in \mathcal{C}$ such that $I\subset C$;
        \item[(b)] for every $C\in\mathcal{C}$, $|C|\leq v(H)-\epsilon L$.
    \end{enumerate}
\end{thm}

Using Theorem~\ref{theorem:container} together with the balanced supersaturation result for expansions of spanning subgraphs of tight trees (Theorem~\ref{theorem:bss_main}), we obtain the following container-type result.

\begin{thm}\label{theorem:container1_main}
    For integers $r\ge k\ge2$, let $\Delta=r-k$. Let $T$ be a tight $k$-tree, let $S$ be a spanning subgraph of $T$ and let $F=S^{(+r)}$. There exist $K,C>0$ such that the following holds for all sufficiently large $n$ and every $t$ with $K\le t\le \binom{n}{r}/n^{r-1}$. Given an $r$-graph $H$ with $n$ vertices and $tn^{r-1}$ edges there exists a collection $\C$ of at most
    $$
    \exp\l(C^{-1}{t^{-\frac{s(F)}{\Delta+1}+1}n^{r-1-\frac{\Delta s(F)}{\Delta+1}}(\log n)^{1+\frac{\Delta s(F)}{\Delta+1}}}\r)
    $$
    subgraphs of $H$ such that
    \begin{enumerate}
        \item[(a)] every $F$-free subgraph of $H$ is a subgraph of some $G\in\C$;
        \item[(b)]$\forall G\in\C$, $|G|\le\l(1-\frac{C}{(\log n)^{\Delta}}\r)tn^{r-1}$.
    \end{enumerate}
\end{thm}

\begin{proof}
    Given that $K$ is sufficiently large, we can apply Theorem~\ref{theorem:bss_main} on $H$. This gives us a constant $C_1>0$ and a collection $\C$ of copies of $F$ in $H$ such that, $\forall 1\le j\le e(F)$,
    $$
    \Delta_j(\C)\le\frac{C_1(\log n)^\Delta|\C|}{tn^{r-1}}A^{-s(F)(j-1)}.
    $$
    where 
    $$
    A=t^{\frac{1}{\Delta+1}}\l(\frac{n}{\log n}\r)^{\frac{\Delta}{\Delta+1}}.
    $$
    We can view $\C$ as a $e(F)$-graph on $E(H)$ whose edges are the edge sets that form the copies of $F$ in $\C$. Then we apply Theorem~\ref{theorem:container} on $\C$ with parameters
    $$
    L=\frac{tn^{r-1}}{C_1(\log n)^{\Delta}}~\text{and}~B=e(H)A^{-s(F)}.
    $$
    This gives us a collection $\C$ of subgraphs of $H$ and a constant $C>0$ such that
    $$
    |\C|\le \exp\l(C^{-1}\log n\cdot e(H)A^{-s(F)}\r)=\exp\l(C^{-1}t^{-\frac{s(F)}{\Delta+1}+1}n^{r-1-\frac{\Delta s(F)}{\Delta+1}}(\log n)^{1+\frac{\Delta s(F)}{\Delta+1}}\r);
    $$
    every $F$-free subgraphs of $H$ is a subgraph of some $G\in\C$; and
    $$
    \forall G\in\C,~e(G)\le \l(1-\frac{C}{(\log n)^{\Delta}}\r)tn^{r-1}.
    $$
\end{proof}

We obtain the following theorem by repeatedly applying Theorem~\ref{theorem:container1_main}.

\begin{thm}\label{theorem:container2_main}
    For integers $r\ge k\ge2$, let $\Delta=r-k$. Let $T$ be a tight $k$-tree, let $S$ be a spanning subgraph of $T$ and let $F=S^{(+r)}$. There exists $K,C>0$ such that for all sufficiently large $n$ and every $t$ with $K\le t\le \binom{n}{r}/n^{r-1}$, there exists a collection $\C$ of $r$-graphs on $[n]$ such that
    \begin{itemize}
        \item[(a)]
        $$
        |\C|\le \exp \l((Ct^{-\frac{s(F)}{\Delta+1}+1}n^{r-1-\frac{\Delta s(F)}{\Delta+1}}(\log n)^{\Delta+2+\frac{\Delta s(F)}{\Delta+1}}\r);
        $$
        \item[(b)] every $F$-free $r$-graph on $[n]$ is a subgraph of some $G\in\C$;
        \item[(c)] $e(G)\le tn^{r-1}$, $\forall G\in\C$.
    \end{itemize}
\end{thm}

\begin{proof}
    Let $H_0=K^r_n$ on $[n]$ and $\C_0=\{H_0\}$. Define $t_0=\binom{n}{r}/n^{r-1}$, $t_i=\exp\l(-\frac{C_1}{(\log n)^{\Delta}}\r)t_{i-1}$ for all $i\ge 1$. Let $m$ be the smallest integer such that $t_m\le t$. It is not hard to check that $m\le O((\log n)^{\Delta+1})$. 
    
    For $0\le i< m$, given $\C_i$ such that $|G|\le t_in^{r-1}$ $\forall G\in\C_i$, we inductively construct $\C_{i+1}$ as follows. Define
    $$
    \C_i^{>}:=\l\{G\in \C_i: |G|> t_{i+1}n^{r-1}\r\},~\text{and~} \C_i^{\le}:=\l\{G\in \C_i: |G|\le t_{i+1}n^{r-1}\r\}.
    $$
    For each $G\in \C_i^{>}$, we apply Theorem~\ref{theorem:container1_main} on $G$ to obtain a collection $\C_G$ of subgraphs of $G$ such that 
    $$
    |\C_G|\le \exp\l(C_1^{-1}{t_{i+1}^{-\frac{s(F)}{\Delta+1}+1}n^{r-1-\frac{\Delta s(F)}{\Delta+1}}(\log n)^{1+\frac{\Delta s(F)}{\Delta+1}}}\r);
    $$
    every $F$-free subgraph of $G$ is a subgraph of some $G'\in\C_G$; and $\forall G'\in\C_G$, 
    $$
    |G'|\le\l(1-\frac{C_1}{(\log n)^{\Delta}}\r)|G|\le \exp\l(-\frac{C_1}{(\log n)^{\Delta}}\r)t_in^{r-1}=t_{i+1}n^{r-1}.
    $$
    Let
    $$
    \C_{i+1}=\C^{\le}_i\cup\l(\bigcup_{G\in \C^>_i}\C_G\r).
    $$
    
    After $m$ repetitions of the above argument, we let $\C=\C_m$. Clearly, every $F$-free subgraph of $K^r_n$ is a subgraph of some $G\in\C$, and $e(G)\le tn^{r-1}$ for every $G\in\C$. Also, for $1\le i\le m$, we have
     $$
     \frac{|\C_{i}|}{|\C_{i-1}|}\le \exp\l(C_1^{-1}{t_{i}^{-\frac{s(F)}{\Delta+1}+1}n^{r-1-\frac{\Delta s(F)}{\Delta+1}}(\log n)^{1+\frac{\Delta s(F)}{\Delta+1}}}\r); 
     $$
    Therefore, for some sufficiently large constant $C>0$,
     $$
     |\C|\le \prod_{i=1}^{m}\frac{|\C_{i}|}{|\C_{i-1}|}\le \exp\l(Ct^{-\frac{s(F)}{\Delta+1}+1}n^{r-1-\frac{\Delta s(F)}{\Delta+1}}(\log n)^{\Delta+2+\frac{\Delta s(F)}{\Delta+1}}\r).
     $$
\end{proof}

\section{Proof of Theorem~\ref{theorem:main}}\label{section:main}
In this section, we prove Theorem~\ref{theorem:main}. For convenience, we restate it here
\begin{thm}
    For integers $r\ge k\ge2$, let $\Delta=r-k$. Let $T$ be a tight $k$-tree, let $S$ be a spanning subgraph of $T$, and let $F=S^{(+r)}$. There exists a constant $C>0$ such that the following holds. Let
    $$
    p_0=Cn^{-s(F)}(\log n)^{\Delta+2+\frac{\Delta s(F)}{\Delta+1}}\text{and}~p_1=n^{-\frac{\Delta}{\Delta+1}s(F)}(\log n)^{\Delta+2+\frac{\Delta s(F)}{\Delta+1}}.
    $$
    Then a.a.s.
    $$
    \ex(G^r_{n,p},F)\le 
    \l\{
    \begin{aligned}
    &Cpn^{r-1},~~~&\text{if}~p\ge p_1;\\
    &Cp^{\frac{s(F)-\Delta-1}{s(F)}}n^{k-1}(\log n)^{\Delta+\frac{(\Delta+2)(\Delta+1)}{s(F)}},~~&\text{if}~p_0\le p\le p_1.
    \end{aligned}
    \r.
    $$
\end{thm}

\begin{proof}
    Let $K=K(F,S,T)$ and $C_1=C_1(F,S,T)>0$ be the constants guaranteed by Theorem~\ref{theorem:container2_main}. Let $X_m$ be the number of $F$-free subgraphs of $G^r_{n,p}$ with $m$ edges. 
    
    If $p\ge p_1$, let
    $$
    m=Cpn^{r-1}.
    $$
    Then we have
    $$
    m\ge Cn^{r-1-\frac{\Delta s(F)}{\Delta+1}}(\log n)^{\Delta+2+\frac{\Delta s(F)}{\Delta+1}}.
    $$
    By Theorem~\ref{theorem:container2_main} with $t=K$, there exists a collection $\C_1$ of subgraphs of $K^r_n$ such that
    $$
    |\C_1|\le \exp\l(C_1K^{-\frac{s(F)}{\Delta+1}+1}n^{r-1-\frac{\Delta s(F)}{\Delta+1}}(\log n)^{\Delta+2+\frac{\Delta s(F)}{\Delta+1}}\r)\le \exp\l(\frac{C_1}{C}K^{-\frac{s(F)-\Delta-1}{\Delta+1}}m\r);
    $$
    $e(G)\le Kn^{r-1}$ for every $G\in\C_1$; and every $F$-free subgraph of $K^r_n$ is a subgraph of some $G\in \C_1$. 
    
    In particular, every $F$-free subgraph of $K^r_n$ with $m$ edges is a subgraph of some $G\in\C_1$. Hence the expectation of $X_m$,
     $$
    \begin{aligned}
    \E[X_m]&\le |\C_1|\binom{Kn^{r-1}}{m}p^m\\
    &\le \exp\l(\frac{C_1}{C}K^{-\frac{s(F)-\Delta-1}{\Delta+1}}m+\log\l(\frac{eKn^{r-1}}{m}\r)\cdot m+\log p\cdot m\r)\\
    &=\exp\l(m\l(\frac{C_1}{C}K^{-\frac{s(F)-\Delta-1}{\Delta+1}}+1+\log K-\log C\r)\r).
    \end{aligned}
    $$
    For $C$ sufficiently large, $\frac{C_1}{C}K^{-\frac{s(F)-\Delta-1}{\Delta+1}}+1+\log K-\log C<0$. Hence $\E[X_m]\to 0$ as $n\to\infty$. By the Markov's Inequality,
    $$
    \P[X_m\ge 1]\le\E[X_m]\to 0,~\text{as}~n\to\infty.
    $$
    This implies that a.a.s.
    $$
    \ex(G^r_{n,p},F)\le m= Cpn^{r-1}.
    $$

    If $p_0\le p\le p_1$, let 
    $$
    m=Cp^{\frac{s(F)-\Delta-1}{s(F)}}n^{k-1}(\log n)^{\Delta+\frac{(\Delta+2)(\Delta+1)}{s(F)}},
    $$ 
    and let 
    $$
    t=Kp^{-\frac{\Delta +1}{s(F)}}n^{-\Delta}(\log n)^{\Delta+\frac{(\Delta+2)(\Delta+1)}{s(F)}}.
    $$ 
    Since $p_0\le p \le p_1$, we have
    $$
    K\le t\le C^{-\frac{\Delta+1}{s(F)}} K n\le \binom{n}{r}/n^{r-1}.
    $$
    By Theorem~\ref{theorem:container2_main}, there exists a collection $\C_2$ of subgraphs of $K^r_n$ such that
    $$
    |\C_2|\le \exp\l(C_1(\log n)^{\Delta+2+\frac{\Delta s(F)}{\Delta+1}}t^{-\frac{s(F)-\Delta-1}{\Delta+1}}n^{r-1-\frac{\Delta s(F)}{\Delta+1}}\r)= \exp\l(\frac{C_1}{C}K^{-\frac{s(F)-\Delta-1}{\Delta+1}}m\r);
    $$
    $e(G)\le tn^{r-1}$ for every $G\in\C_2$; and every $F$-free subgraph of $K^r_n$ is a subgraph of some $G\in \C_2$. 
    
    In particular, every $F$-free subgraph of $K^r_n$ with $m$ edges is a subgraph of some $G\in\C_2$. Hence the expectation of $X_m$,
    $$
    \begin{aligned}
    \E[X_m]&\le |\C_2|\binom{tn^{r-1}}{m}p^m\\
    &\le \exp\l(\frac{C_1}{C}K^{-\frac{s(F)-\Delta-1}{\Delta+1}}m+\log\l(\frac{etn^{r-1}}{m}\r)\cdot m+\log p\cdot m\r)\\
    &=\exp\l(m\l(\frac{C_1}{C}K^{-\frac{s(F)-\Delta-1}{\Delta+1}}+1+\log K-\log C\r)\r)\to0~\text{as}~n\to\infty.
    \end{aligned}
    $$
    By the Markov's Inequality,
    $$
    \P[X_m\ge 1]\le\E[X_m]\to 0,~\text{as}~n\to\infty.
    $$
    This implies that a.a.s.
    $$
    \ex(G^r_{n,p},F)\le m=Cp^{\frac{s(F)-\Delta-1}{s(F)}}n^{k-1}(\log n)^{\Delta+\frac{(\Delta+2)(\Delta+1)}{s(F)}}.
    $$ 
\end{proof}

\section{Applications of Theorem~\ref{theorem:main} and tight constructions}\label{section:construction}
We first introduce two simple constructions.
\begin{prop}\label{proposition:random}
    Let $H$ be an $r$-graph. If $\omega(n^{-r})\le p\le o(n^{-s(H)})$, then a.a.s.
    $$
    \ex(G^r_{n,p},H)=\l(1+o(1)\r)p\binom{n}{r}.
    $$

    If $p\ge n^{-s(H)}$, then a.a.s.
    $$
    \ex(G^r_{n,p},H)\ge \Omega\l(n^{r-s(H)}\r).
    $$
\end{prop}

\begin{proof}[Sketch proof]
     Since $p\ge \omega(n^{-r})$, by the Chernoff bound, a.a.s. $e(G^r_{n,p})=(1+o(1))p\binom{n}{r}$. Hence
     $$
     \ex(G^r_{n,p},H)\le \l(1+o(1)\r)p\binom{n}{r}.
     $$

    Let $G'$ be a subgraph of $G^r_{n,p}$ obtained by arbitrarily deleting an edge for each copy of $H$ in $G^r_{n,p}$, and let $N(H,G^r_{n,p})$ be the number of copies of $H$ in $G^r_{n,p}$. Since $p\le o(n^{-s(H)})$, it is not hard to check that
    $$
    \E[N(H,G^r_{n,p})]=o\l(pn^r\r).
    $$
    Hence, by the Markov inequality,
    $$
    \P[N(H,G^r_{n,p}\ge o(pn^r)]=o(1).
    $$
    Therefore, a.a.s.
    $$
    \ex(G^r_{n,p},H)\ge e(G')\ge e(G^r_{n,p})-N(H,G^r_{n,p})=\l(1+o(1)\r)p\binom{n}{r}.
    $$

    If $p\ge n^{-s(H)}$, since $\ex(G^r_{n,p},H)$ is non-decreasing with respect to $p$, let $c>0$ be a sufficiently small constant, we have a.a.s.
    $$
    \ex(G^r_{n,p},H)\ge \ex(G^r_{n,cn^{-s(H)}},H)\ge \Omega\l(n^{r-s(H)}\r).
    $$
\end{proof}

\begin{prop}\label{proposition:star}
    Let $H$ be an $r$-graph such that $\cap_{e\in E(H)}e=\emptyset$. If $p=\omega(n^{-r+1})$, then a.a.s.
    $$
    \ex(G^r_{n,p},H)\ge \Omega(pn^{r-1}).
    $$
\end{prop}

\begin{proof}[Sketch proof]
    Fix a vertex $v$ of $G^r_{n,p}$. Let $G'$ be a subgraph of $G^r_{n,p}$ whose edge set consists of all edges of $G^r_{n,p}$ containing $v$. Then by the Chernoff bound, a.a.s.
    $$
    \ex(G^r_{n,p},H)\ge e(G')\ge\Omega(pn^{r-1}).
    $$
\end{proof}

By Theorem~\ref{theorem:main}, Proposition~\ref{proposition:random} and Proposition~\ref{proposition:star}, we immediately obtain essentially tight bounds for the random Tur\'an numbers of expansions of tight tree, that is, Theorem~\ref{theorem:main_tighttree}.

\begin{proof}[Proof of Theorem~\ref{theorem:main_tighttree}]
By definition, $s(T^{(+r)})=\Delta+1$. Hence by Theorem~\ref{theorem:main} and Proposition~\ref{proposition:random}, we have the upper bound in Theorem~\ref{theorem:main_tighttree}. Moreover, by Proposition~\ref{proposition:random} and Proposition~\ref{proposition:star}, we know that the upper bound is tight.
\end{proof}

For almost all degenerate graphs/hypergraphs $H$ that have been studied, $\ex(G^r_{n,p},H)$ has a ``flat'' middle range. Surprisingly, we show that for all $r\ge k\ge 3$, $\ex(G^r_{n,p},K_k^{k-1(+r)})$ has a ``non-flat'' middle range. Upper bound comes from Theorem~\ref{theorem:main}. Now we introduce a construction which gives a tight lower bound. We make use of the following constructions introduced by Gowers and Janzer~\cite{gowers2021generalizations}, which generalized the famous construction of Ruzsa and Szemerédi~\cite{ruzsa1978}.
\begin{thm}[\cite{gowers2021generalizations,ruzsa1978}]\label{theorem:GJ}
For $r>k\ge 2$ and $n\ge 1$, there exists a graph on $n$ vertices with the following two properties:
\begin{itemize}
    \item[(i)] It has $n^ke^{-O(\sqrt{\log n})}$ subgraphs isomorphic to $K_r$;
    \item[(ii)] For any $t$ with $k<t\le r$ and any subgraph $G_1$ isomorphic to $K_k$, if there exist a subgraph $G_2$ isomorphic to $K_t$ and a subgraph $G_3$ isomorphic to $K_r$ such that $G_1\subset G_2$ and $G_1\subset G_3$, then $G_2\subset G_3$.
\end{itemize}
\end{thm}
Note that property $(ii)$ is slightly stronger than the original Theorem 1.2 in~\cite{gowers2021generalizations} which states ``every $K_k$ is contained in at most one $K_r$". In fact, property $(ii)$ is inherently implied by the proof of Lemma 3.1 in~\cite{gowers2021generalizations}. 

Let $G_{m,r,k}$ be the graph on $m$ vertices guaranteed by Theorem~\ref{theorem:GJ}. Consider an $r$-graph $H_{m,r,k}$ on $V(G_{m,r,k})$ whose edges are the vertex sets of copies of $K_r$ in $G_{m,r,k}$.

\begin{prop}\label{proposition:GJ1}
For $r>k\ge 2$ and $m\ge1$, $H_{m,r,k}$ has the following properties:
\begin{itemize}
    \item[(i)] $e(H_{m,r,k})\ge m^k\exp\l(-O(\sqrt{\log m})\r)$;
    \item[(ii)] Any two edges intersect in at most $k-1$ vertices;
    \item[(iii)] It does not contain any subgraph isomorphic to $K_{k+1}^{k(+r)}$.
\end{itemize}
\end{prop}
\begin{proof}
    Property (i) is trivial. Suppose there are two edges intersect in $k$ vertices, then there exists a copy of $K_k$ contained in two different copies of $K_r$ in $G_{m,r,k}$. This contradicts the property of $G_{m,r,k}$. Similarly, if $H_{m,r,k}$ contains a copy of $K_{k+1}^{k(+r)}$, then in $G_{m,r,k}$ there is a copy of $K_k$ contained in a copy of $K_{k+1}$ and a copy of $K_r$, but the copy of $K_{k+1}$ is not contained in the copy of $K_r$, contradicting the property of $G_{m,r,k}$.
\end{proof}

Let $H_{m,r,k}(n)$ be a hypergraph on $n$ vertices obtained from $H_{m,r,k}$ by replacing each vertex with an independent set of size $\lceil n/m\rceil$ or $\lfloor n/m \rfloor$; and replacing each edge with a complete $r$-partite $r$-graph. More precisely, let $\{1,\dots,m\}$ be the vertices of $H_{m,k,r}$ and let $V_1\cup\dots\cup V_m$ be the vertex set of $H_{m,r,k}(n)$. Define a function $f:V_1\cup\dots\cup V_m\to \{1,\dots,m\}$ such that $f(v)=w$ if and only if $v\in V_w$. Then $e=\{v_1,\dots, v_r\}$ is an edge of $H_{m,r,k}(n)$ if and only if $f(e)=\{f(v_1),\dots, f(v_r)\}$ is an edge of $H_{m,k,r}$. We call $H_{m,r,k}(n)$ a {\em balanced blowup} of $H_{m,r,k}$ on $n$ vertices.

\begin{prop}\label{proposition:GJ2}
If $H_{m,r,k}(n)$ contains a subgraph $G$ isomorphic to $K^{k(+r)}_{k+1}$, then $\forall e_1,e_2\in E(G)$, $f(e_1)=f(e_2)$.
\end{prop}
\begin{proof}
    Labeled the edges of $G$ as $e_1,\dots,e_{k+1}$ and let $e_i=\{v_1,\dots,v_{k+1},w_{i,1},\dots,w_{i,r-k}\}\setminus\{v_i\}$, $1\le i\le k+1$. Without lost of generality, we suppose for contradiction that $f(e_1)\not=f(e_2)$, then, by property $(ii)$ in Proposition~\ref{proposition:GJ1}, $|f(e_1)\cap f(e_2)|\le k-1$. Since $|e_1\cap e_2|=k-1$, we have $|f(e_1)\cap f(e_2)|\ge k-1$. Hence $|f(e_1)\cap f(e_2)|= k-1$. This implies that $f(v_1)\not\in f(e_1)$. Note that $v_1\in e_i$, and hence $f(v_1)\in f(e_i)$, $\forall i\not=1$. This implies $f(e_1)\not=f(e_i)$, and hence, $|f(e_1)\cap f(e_i)|=k-1$, $\forall i\not=1$. By symmetry, $|f(e_i)\cap f(e_j)|=k-1$, $\forall i\not=j$. This means that $f(e_1),\dots,f(e_{k+1})$ form a copy of $K^{k(+r)}_{k+1}$ in $H_{m,r,k}$, which contradicts property $(iii)$ in Proposition~\ref{proposition:GJ1}.
\end{proof}

\begin{thm}\label{theorem:simplex_lower}
    For $r\ge k\ge 3$, let $\Delta=r-k$, there exists $c>0$ such that if $p\ge n^{-\Delta-\frac{k}{k-1}}\exp(c\sqrt{\log n})$, then a.a.s.
    $$
    \ex(G^r_{n,p}, K^{k-1(+r)}_k)\ge p^{\frac{1}{(\Delta+1) (k-1)+1}}n^{k-1}\exp(-c\sqrt{\log n}).
    $$
\end{thm}

\begin{proof}
    Let
    $$
    m=p^{\frac{k-1}{(\Delta+1)(k-1)+1}}n\exp(\sqrt{\log n}),
    $$
    and let $H=H_{m,r,k-1}(n)$.
    Let $X$ be the number of edges in $H\cap G^r_{n,p}$ and let $Y$ be the number of copies of $K^{k-1(+r)}_k$ in $H\cap G^r_{n,p}$. By Proposition~\ref{proposition:GJ1}, 
    $$
    \E[X]=\Theta\l(e(H_{m,r,k-1})\l(\frac{n}{m}\r)^rp\r)=\Theta\l(\frac{e(H_{m,r,k-1})}{m^{k-1}}p^{\frac{1}{(\Delta+1)(k-1)+1}}n^{k-1}\exp\l(-(\Delta+1)\sqrt{\log n}\r)\r).
    $$
    It is not hard to check that $\E[X]\rightarrow\infty$ as $n\rightarrow\infty$ given that $c$ is sufficiently large. Hence, by the Chernoff bound, a.a.s. $X\ge \E[X]/2$.
    
    On the other hand, by Proposition~\ref{proposition:GJ2},
    $$
    \E[Y]=\Theta\l(e(H_{m,r,k-1})\l(\frac{n}{m}\r)^{k(\Delta+2)}p^{k}\r)=\Theta\l(\frac{e(H_{m,r,k-1})}{m^{k-1}}p^{\frac{1}{(\Delta+1)(k-1)+1}}n^{k-1}\exp\l(-(\Delta k+k+1)\sqrt{\log n}\r)\r).
    $$
    Note that $\E[Y]=o(\E[X])$. Hence, by the Markov's inequality,
    $$\P[Y\ge \E[X]/4]\le o(1),$$
    which implies a.a.s. $Y\le \E[X]/4$.
    Let $H'$ be an $K^{k-1(+r)}_k$-free subgraph of $H$ obtained by deleting an edge in each copy of $K^{k-1(+r)}_k$. Then a.a.s.
    $$
    \ex(G^r_{n,p},K^{k-1(+r)}_k)\ge e(H')\ge X-Y\ge \frac{\E[X]}{4}\ge {p^{\frac{1}{(\Delta+1)(k-1)+1}}}n^{k-1}\exp\l(-O\l(\sqrt{\log n}\r)\r).
    $$
\end{proof}

\begin{proof}[Proof of Theorem~\ref{Theorem:main_clique}]
    Observe that $K^{k-1(+k)}_k$ is contained in a tight $k$-tree as spanning subgraph: begin with a $k$-edge $e$ and then expand each $(k-1)$-shadow of $e$ into a $k$-edge; this gives a tight $k$-tree containing $K^{k-1(+k)}_k$ as spanning subgraph. It is not hard to check that $s(K^{k-1(+r)}_k)=\Delta+\frac{k}{k-1}$. Hence we can apply Theorem~\ref{theorem:main} on $K^{k-1(+r)}_k$, which implies the desired upper bound. Lower bound comes from Proposition~\ref{proposition:star} and Theorem~\ref{theorem:simplex_lower}.
\end{proof}

\section*{Acknowledgement}
The author would like to express gratitude to Jacques Verstraëte for introducing him to the random Turán problems. Additionally, the author would like to thank Hehui Wu for helpful discussion in the early stage of this paper. Last but not least, the author would like to thank Barnabas Janzer for confirming the correctness of property (ii) in Theorem~\ref{theorem:GJ}.

\bibliographystyle{abbrv}
\bibliography{refs}

\end{document}